\documentclass[reqno,12pt]{amsart}

\usepackage[T1]{fontenc}
\usepackage{lmodern}

\usepackage[numbers,sort&compress]{natbib}
\usepackage{graphicx}
\usepackage{hyperref}
\usepackage{amsmath,amssymb,amsthm}
\usepackage{amsaddr}
\usepackage{mathtools}
\usepackage{textcomp}
\usepackage{algorithm}
\usepackage{algpseudocode}
\usepackage{array}
\usepackage[labelfont={rm}, justification=centering]{subfig}
\usepackage{booktabs}

\usepackage[svgnames]{xcolor}

\calclayout

\DeclareMathAlphabet{\mathbsl}{OT1}{cmr}{bx}{sl}  

\newcommand{\M}[1]{\mathbsl{#1}}

\newcommand{\MS}[1]{\boldsymbol{#1}}

\newcommand{\D}    {\operatorname{d}}
\newcommand{\sgn}  {\operatorname{sgn}}

\newcommand{\transpose}[1]{#1^t}

\newcommand{\nel}       {n_{\textsc{e}}}
\newcommand{\ns}  [1][] {n_{\textsc{s}#1}}
\newcommand{\no}  [1][] {n_{\textsc{o}#1}}
\newcommand{\np}        {n_{\textsc{p}}}
\newcommand{\deltaO}    {\delta_{\textsc{o}}}
\newcommand{\nit}       {n_{\textsc{i}}}

\newcommand{\cs}        {c_{\textsc{s}}}

\newcommand{\IOP}{\boldsymbol{\mathcal{I}}}

\newcommand{\tsup}[1]{$^{#1}$}

\begin{document}

\title%
[Schwarz smoothers for spectral element multigrid]%
{Nonuniformly weighted Schwarz smoothers for spectral element multigrid}

\author{J\"org Stiller}
\address{%
  Technische Universit\"at Dresden, Institute of Fluid Mechanics and
  Center for Advancing Electronics Dresden, 01062 Dresden, Germany}
\email{joerg.stiller@tu-dresden.de}

\thanks{%
Article accepted for publication in J. Sci. Comput.
}

\begin{abstract}
A hybrid Schwarz/multigrid method for spectral element solvers to the Poisson equation in $\mathbb R^2$ is presented.
It extends the additive Schwarz method studied by J. Lottes and P. Fischer (J. Sci. Comput. 24:45--78, 2005) by introducing nonuniform weight distributions based on the smoothed sign function.
Using a V-cycle with only one pre-smoothing, the new method attains logarithmic convergence rates in the range from 1.2 to 1.9, which corresponds to residual reductions of almost two orders of magnitude.
Compared to the original method, it reduces the iteration count by a factor  of 1.5 to 3, leading to runtime savings of about 50 percent.
In numerical experiments the method proved robust with respect to the mesh size and polynomial orders up to 32.
Used as a preconditioner for the (inexact) CG method it is also suited for anisotropic meshes and easily extended to diffusion problems with variable coefficients.
\end{abstract}

\keywords{%
Multigrid method,
Schwarz methods,
spectral element method,
$p$-version finite element method.}

\maketitle


\section{Introduction}
\label{sec:intro}

High-order finite element methods (FEM) enjoy an increasing interest in computational science and engineering. They include $hp$-FEM, spectral element methods (SEM) as well as discontinuous Galerkin methods \cite{KS05,DFM02}. The motive that drives the development of high-order methods lies in their potential to deliver accuracy with lower cost in comparison to the first and second order methods used in common simulation tools \cite{Wang2013}. However, realizing this advantage in practice is a formidable task. Along with curvilinear mesh generation, the provision of efficient solvers for the resulting algebraic equation systems remains the main challenge.

Projection methods for incompressible flow, or implicit discretization of diffusion terms lead to a sequence of linear elliptic problems which are related or equivalent to the Poisson equation or, more generally, the Helmholtz equation \cite{Guermond2006}.
Fast solvers for such equations are therefore a crucial ingredient of competitive high-order methods and, hence, have been in focus of research for almost 30 years
\cite{RP87,Hei88,MM88,FL04,Kan04,BZ05,LF05,Ols07,KT08,LNS08,PR09,Mit10,JK11,BBS12,HSN13}.
For Helmholtz or Poisson problems discretized on regular meshes, efficient multigrid (MG) techniques have been developed recently.
\citet{LF05} proposed additive Schwarz smoothers based on extended element domains, which attain residual reductions of approximately 0.2 within one sweep.
They found that weighting the overlapping Schwarz updates by the inverse of the counting matrix, which corresponds to taking the arithmetic mean, plays a crucial role in obtaining multigrid-like iteration counts.
A detailed analysis of the method was given in \cite{LNS08}.
\citet{JK11} presented a similar multigrid approach for the $p$-finite element method on locally refined Cartesian meshes.
They used a multiplicative Schwarz smoother on element domains which possess only a minimal overlap confined to the element boundaries.
\citet{HSN13} developed a $p$-multigrid method based on static condensation which, apart from pre- and post-processing, reaches linear complexity.
The proposed block smoother can be classified as an additive Schwarz method using a monotonic increasing shape function for blending the overlapping updates.
Using this smoother the multigrid method attained convergence rates of about 0.02 combined with a run-time efficiency that comes close to fast direct finite difference solvers.
The success of this approach inspired us to extend the idea of nonuniform weighting to the full, "uncondensed" problem and thus led to the present work.
The primary goal is to show how nonuniform weighting can be used to boost the performance of high-order spectral-element multigrid techniques.
Further, we investigate the influence of the overlap width, smoothing strategies, additive versus multiplicative Schwarz methods and Krylov acceleration on robustness and efficiency.
In addition to this, we consider the extension to diffusion problems with variable coefficients.

The remainder of the paper is organized as follows:
Section~\ref{sec:discretization} provides a brief description of the spectral element discretization.
Section~\ref{sec:methods} presents the solution techniques, namely the weighted additive and multiplicative Schwarz methods, the $p$-multigrid method and the inexact multigrid-preconditioned conjugate gradient method.
Section~\ref{sec:results} proceeds with the discussion of numerical experiments for assessing the solution methods and application to variable diffusion.
Finally, Section~\ref{sec:conclusions} concludes the paper.


\section{Discretization}
\label{sec:discretization}

As the model problem we consider the Poisson equation
\begin{equation}
\label{eq:poisson}
-\nabla^2 u = f
\end{equation}
in the rectangular domain \mbox{$\Omega=[0,\ell_x]\times[0,\ell_y]$} with periodic boundaries.
For discretization $\Omega$ is decomposed into $\nel=n_x \times n_y$ rectangular elements
$\Omega^{mn}$
with dimensions
\mbox{$\Delta x = \ell_x/n_x$} and
\mbox{$\Delta y = \ell_y/n_y$}. 
%
In each element the solution is approximated as
\begin{equation}
\label{eq:ansatz}
u(x,y)|_{\Omega^{mn}} \simeq
\sum_{i,j=0}^{p} u^{mn}_{ij} \, \varphi_i\big(\xi^m(x)\big) \, \varphi_j\big(\eta^n(y)\big)
\end{equation}
where
$\varphi_i$ are the Lagrange polynomials to the Gauss-Lobatto-Lengendre (GLL) points
$\{\xi_i\}_{i=0}^{p}$
in the one-dimensional standard region \mbox{$\hat\Omega = [-1,1]$}
and
$\xi^m(x)$,
$\eta^n(y)$
the mapping of coordinates from $\Omega^{mn}$ to $\hat\Omega$
\cite{KS05,DFM02}.
Concatenation of the element coefficients
\mbox{$\M u^{mn} = [u^{mn}_{ij}]$}
and enforcing continuity for shared vertices and edges
yields the unique global coefficients $\M u$ \cite[see, e.g.][pp.\ 191--194]{DFM02}.
Application of the Galerkin spectral element method leads to the discrete equations
\begin{equation}
\label{eq:discrete}
\M A \M u = \M f \,.
\end{equation}
As a consequence of the tensor product ansatz \eqref{eq:ansatz} and the Cartesian mesh,
the global system matrix in Eq.~\eqref{eq:discrete} assumes the tensor product form
\begin{equation}
\label{eq:system matrix}
\M A = \M M_y \otimes \M L_x + \M L_y \otimes \M M_x \, ,
\end{equation}
where
$\M M_\ast$ and
$\M L_\ast$
represent the one-dimensional mass and stiffness matrices for directions
\mbox{$\ast=x,y$}, respectively.
The detailed structure of these operators and underlying spectral element techniques are well described in literature
\cite{DFM02,KS05}
and therefore deliberately skipped here.


\section{Solution methods}
\label{sec:methods}

For solving Eq.~\eqref{eq:discrete} we consider polynomial multigrid (MG) and
multigrid-preconditioned conjugate gradients (MGCG).
Both approaches rely on Schwarz methods for smoothing.
We first present the Schwarz methods and then sketch MG and MGCG.

\subsection{Schwarz methods}
\label{sec:Schwarz}

\begin{figure}
\centering
\includegraphics[width=0.65\textwidth]{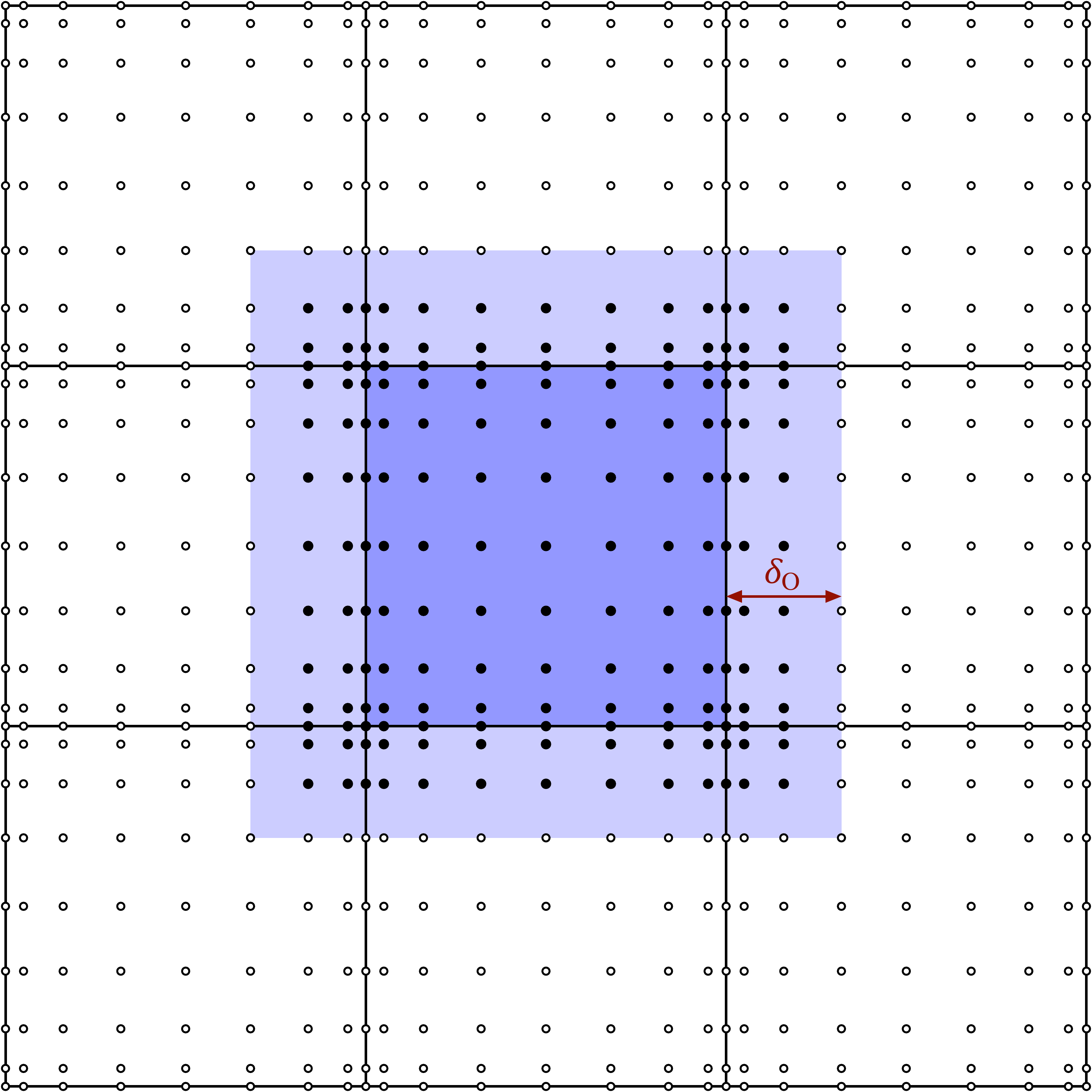}
\caption{Example of a subdomain used with the Schwarz method. The shaded area represents $\Omega_s$ and the dark region in its center the corresponding element. The circles are the GLL nodes for polynomial order $p=8$. Filled circles indicate the nodes that are solved for and updated; $\deltaO$ is the overlap width.}
\label{fig:subdomain}
\end{figure}

Schwarz methods are iterative domain decomposition techniques which improve the approximate solution by parallel or sequential subdomain solves, leading to additive or multiplicative methods, respectively.
Following \cite{LF05} we use extended element regions as the subdomains.
Figure~\ref{fig:subdomain} illustrates how the subdomain $\Omega_{s}$ results from the corresponding element domain $\Omega^{mn}$ by attaching a rectangular strip matching the overlap width $\deltaO$.
As consequence, $\Omega_{s}$ adopts $\no$ layers of additional nodes from the neighbor elements.
Note, however, that we exclude the outer layer of nodes located on $\partial\Omega_{s}$.
For definiteness we define the overlap width in terms of $\no$ and the GLL points,
\begin{equation}
\deltaO = \xi_{\no + 1} + 1 \, .
\end{equation}

To derive a local correction to some approximate solution $\M{\tilde u}$ we first convert Eq.~\eqref{eq:discrete} into the equivalent residual form
\begin{equation}
\label{eq:residual form}
\M A \Delta \M u = \M f - \M A \M{\tilde u} = \M{\tilde r} \, ,
\end{equation}
where
${\Delta \M u = \M u - \M {\tilde u}}$.
Further we introduce the restriction operator $\M R_s$ such that
\mbox{$\M u_s = \M R_s \M u$} gives the coefficients associated with $\Omega_{s}$.
Conversely, the transposed restriction operator globalizes any local coefficients by adding zeros for exterior nodes.
With these prerequisites the correction contributed by $\Omega_{s}$ is defined as the solution of the subproblem
\begin{equation}
\label{eq:subproblem}
\M A_{ss} \Delta \M u_s = \M r_s \, ,
\end{equation}
where
\mbox{$\M A_{ss} = \M R_s \M A \transpose{\M R}_s$}
represents the restricted system matrix and
\mbox{$\M r_s = \M R_s \M{\tilde r}$}
the restricted residual.
Due to the rectangular shape of the subdomain, $\M A_{ss}$ inherits the tensor product structure of $\M A$.
Using the fast diagonalization technique developed by \citet{Lyn64} and adopted for SEM e.g. in \cite{CD95}, the inverse subdomain operator can be expressed in the form
\begin{equation}
\M A_{ss}^{-1}
= (\M S_y \otimes \M S_x)
  (\M I \otimes \MS \Lambda_x + \MS \Lambda_y \otimes \M I)^{-1}
  (\transpose{\M S_y} \otimes \transpose{\M S_x}),
\end{equation}
where $\M I$ is the unity matrix,
$\M S_{\ast}$ the matrix of eigenvectors to the generalized eigenproblem for the restricted one-dimensional stiffness and mass matrices,
and $\MS \Lambda_{\ast}$ the diagonal matrix of eigenvalues
for directions $\ast=x,y$.
With equidistant meshes, as in the present case, the operators are identical for all subdomains and, hence, the cost for their pre-computation becomes negligible.
Exploiting the tensor-product structure of the inverse, the solution to a single subdomain,
\mbox{$\Delta\M u_s = \M A_{ss}^{-1} \M r_s$},
can be evaluated with just $\Theta\big(2(p+1+2\no)^3\big)$ operations.

There exist several options for combining the local solutions.
We consider a weighted version of the additive Schwarz method and the multiplicative Schwarz method.
The multiplicative Schwarz method solves the subproblems \eqref{eq:subproblem} consecutively while continually updating the residual.
Note that, in general, one multiplicative Schwarz iteration corresponds to the application of a non-symmetric linear operator, albeit $\M A$ is symmetric.
However, for an even number of steps, the method is symmetrized by reversing the order of subdomains in each step, which leads to Algorithm~\ref{alg:schwarz:mult}.

\begin{algorithm}[t]
\caption{Multiplicative Schwarz method}
\label{alg:schwarz:mult}
\begin{algorithmic}[1]
\Function{MSchwarz}{$\M u, \M f, \nit$}
   \For{$i=1,\nit$}
      \For{$e=1,\nel$}
      \State \( s \gets
                \left\{
                   \begin{array}{l@{\quad}l}
                      e             & i \ \text{odd}   \\
                      \nel + 1 - e  & i \ \text{even}
                   \end{array}
                \right.
             \)
      \State $\M r \gets \M f - \M A \M u$
      \State $\M u \gets \M u + \transpose{\M R}_s \M A_{ss}^{-1} \M R_s \M r$
      \EndFor
   \EndFor
   \State \textbf{return} $\M u$
\EndFunction
\end{algorithmic}
\end{algorithm}

\begin{algorithm}[t]
\caption{Weighted additive Schwarz method}
\label{alg:schwarz:wadd}
\begin{algorithmic}[1]
\Function{WSchwarz}{$\M u, \M f, \nit$}
   \For{$i=1,\nit$}
      \State $\M r \gets \M f - \M A \M u$
      \State $\M u \gets \M u
              + \sum_{s=1}^{\nel}
                   \transpose{\M R}_s \M W \M A_{ss}^{-1} \M R_s \M r $
   \EndFor
   \State \textbf{return} $\M u$
\EndFunction
\end{algorithmic}
\end{algorithm}

\begin{figure}
\centering
\includegraphics[width=0.65\textwidth]{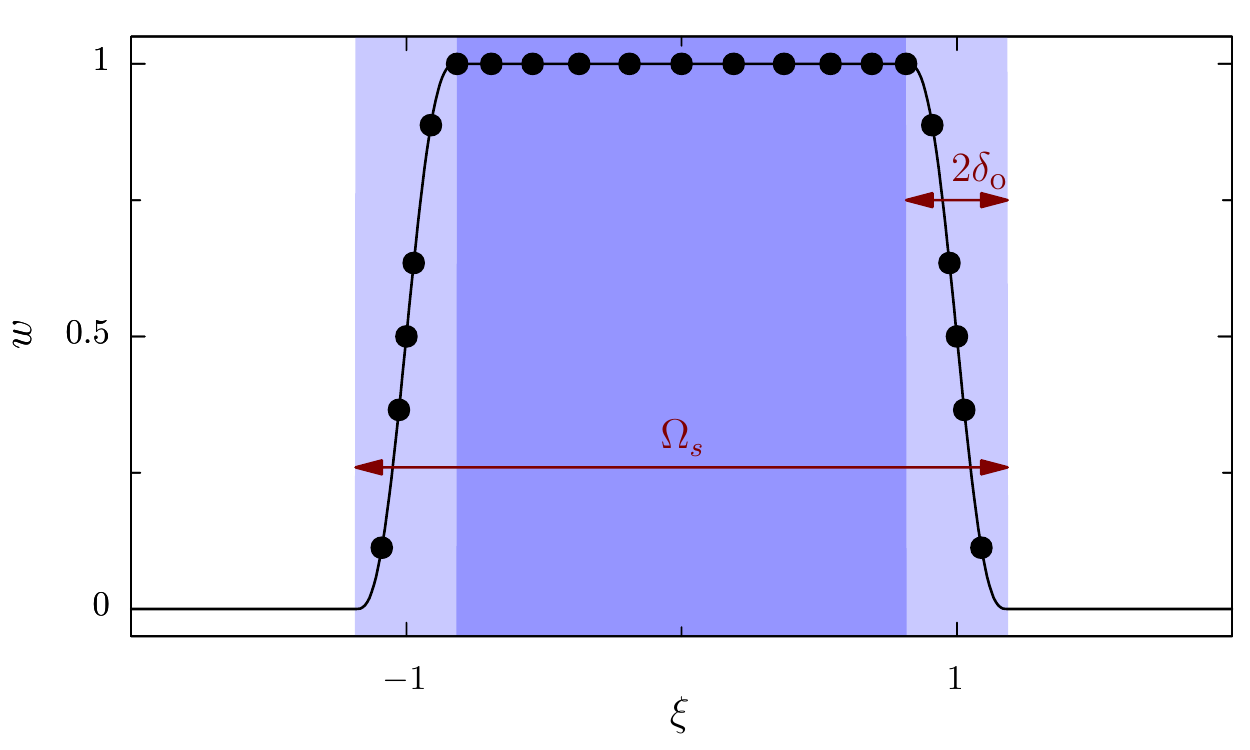}
\caption{One-dimensional weight distribution for elements of order $p=16$ with an overlap of $\no=2$ points using a quintic shape function.
The core region and the overlap zone of the subdomain are shaded in dark and light blue, respectively.
Filled circles indicate the node positions.
}
\label{fig:weight_distribution}
\end{figure}

The weighted additive Schwarz method determines all local corrections independently and computes the global correction as a linear combination of these results, i.e.
\begin{equation}
\label{eq:correction:global}
\Delta \M u \simeq \sum_{s} \transpose{\M R_{s}} (\M W \Delta \M u_{s}) \, ,
\end{equation}
where $\M W$ is a diagonal local weight matrix.
Application of Eq.~\eqref{eq:correction:global} leads to Algorithm~\ref{alg:schwarz:wadd}.
Note that \mbox{$\M W = \M I$} recovers the classical additive Schwarz method.
The arithmetic mean employed in \cite{LF05} is obtained by choosing \mbox{$\M W = \M R_s \M C^{+} \transpose{\M R}_s$}, where $\M C^{+}$ is the pseudoinverse of the counting matrix $\M C = \sum_s \transpose{\M R}_s \M R_s$.
We propose a more flexible approach which elevates the weights gradually from zero at the border to one in the core zone.
Due to the regular shape of $\Omega_s$ the weights can be cast in the tensor product form ${\M W = \M W_y \otimes \M W_x}$.
The one-dimensional weight distributions $\M W_{\ast}$ are generated from the continuous weighting function
\begin{equation}
\label{eq:weight function}
w_{\kappa}(\xi) = \frac{1}{2} \left[
                            \phi_{\kappa}\left(\frac{\xi + 1}{\deltaO}\right)
                          - \phi_{\kappa}\left(\frac{\xi - 1}{\deltaO}\right)
                     \right]
\, ,
\end{equation}
where
$\xi$ is the 1D standard coordinate extended beyond $\hat \Omega$ and
$\phi_{\kappa}$ is a weakly monotonic increasing shape function.
In particular we consider the shape functions $\phi_{\kappa}$ with ${\kappa}\in\{1,3,5,\dots\}$ defined as
\begin{equation}
\label{eq:shape function}
\phi_{\kappa}(x)
= \left\{
  \begin{array}{l@{\qquad}l}
    \hat\phi_{\kappa}(x) &  x \in \hat\Omega \\
    \sgn(x)                 &  \text{else}
  \end{array}
  \right.
\end{equation}
where $\hat\phi_{\kappa}$ is a polynomial of degree $i$ satisfying the conditions
\begin{subequations}
\begin{align}
& \hat\phi_{\kappa}(\pm 1) = \pm 1  \\
& \frac{\D^k \hat\phi_{\kappa}}{\D x^k}(\pm 1) = 0, \quad 0 < k \le ({\kappa}-1)/2
\, .
\end{align}
\end{subequations}
The $\hat\phi_{\kappa}$ are strictly monotonic in $(-1,1)$ and generate a smooth transition of the weight function in the overlap zone, as exemplified in Fig.~\ref{fig:weight_distribution} for the quintic case.
By increasing the polynomial degree the shape function converges toward the sign function, which translates into a top hat weighting function.
We remark that omitting the shape function in Eq.~\eqref{eq:weight function} yields the arithmetic mean.
For reference, Table~\ref{tab:weighting} summarizes all weight functions used in the numerical experiments.

\begin{table}
\renewcommand{\arraystretch}{1.3}
\caption{Weight functions (WF) and related shape functions}
\label{tab:weighting}
\centering
\begin{tabular}{cll}
\hline\noalign{\smallskip}
WF & shape function                                  & method  \\
\noalign{\smallskip}\hline\noalign{\smallskip}
$w_1$  & $\hat\phi_{1} = x$                                 & linear  \\
$w_3$  & $\hat\phi_{3} = (3x - x^3) / 2$                    & cubic   \\
$w_5$  & $\hat\phi_{5} = (15x - 10x^3 +  3x^5) / 8$         & quintic \\
$w_7$  & $\hat\phi_{7} = (35x - 35x^3 + 21x^5 - 5x^7) / 16$ & 7\tsup{th} order \\
$w_t$  & $\hat\phi_{t} = \sgn(x)$                           & top hat \\
$w_a$  & $\hat\phi_{a} = 0$                                 & arithmetic mean \\
\noalign{\smallskip}\hline
\end{tabular}
\end{table}

\subsection{Multigrid}
\label{sec:mg}

For MG we define a series of polynomial levels $\{p_l\}$
with
\mbox{$p_l = 2^l$}
increasing from $1$ at \mbox{$l=0$} to $p$ at top level $L$.
Correspondingly,
$\M u_l$ denotes the global coefficients and
$\M A_l$ the system matrix on level $l$.
On the top level we have
\mbox{$\M u_L = \M u$} and
\mbox{$\M A_L = \M A$},
whereas on lower levels $\M u_l$ is the defect correction and $\M A_l$ the counterpart of $\M A$ obtained with elements of order $p_l$.
For transferring the correction from level \mbox{$l-1$} to level $l$ we use the embedded interpolation operator $\IOP_l$, and for restricting the residual its transpose.
These ingredients allow to build a multigrid V-cycle, which is summarized in Algorithm~\ref{alg:V-cycle}.
Both, the multiplicative and the weighted additive Schwarz method stated in Algorithm~\ref{alg:schwarz:mult} and \ref{alg:schwarz:wadd}, respectively, can serve as the \textsc{Smoother}.
The number of pre- and post-smoothing steps, $\ns[1,l]$ and $\ns[2,l]$, can differ from level to level to allow variable V cycles \cite{Bra93}.
Line~\ref{alg:V-cycle:coarse} of Algorithm~\ref{alg:V-cycle} defines the coarse grid solution by means of the pseudoinverse $\M A_0^{+}$.
In our implementation the coarse problem is solved using the conjugate gradient method. To achieve convergence in spite of singularity, the right side is projected to the null space of $\M A_0$, as proposed by \citet{Kaa88}.

\begin{algorithm}[ht]
\caption{Multigrid V-cycle.}
\label{alg:V-cycle}
\begin{algorithmic}[1]
\Function{MultigridCycle}{$\M u$, $\M f$, $\M\ns$}
   \State $\M u_L \gets \M u$
   \State $\M f_L \gets \M f$
   \For{$l=L,1$ \textbf{step} $-1$}
       \State $\M u_l \gets$ \Call{Smoother}{$\M u_l$, $\M f_l$, $\ns[1,l]$}
         \Comment{Pre-smoothing}
       \State $\M f_{l-1} \gets \transpose{\IOP}_l (\M f_l - \M A_l \M u_l)$
         \Comment{Residual restriction}
   \EndFor
   \State $\M u_0 \gets \M A_0^{+} \M f_0$
     \label{alg:V-cycle:coarse}
     \Comment{Coarse grid solution}
   \For{$l=1,L$}
       \State $\M u_l \gets \M u_l + \IOP_l \M u_{l-1}$
         \Comment{Correction prolongation}
       \State $\M u_l \gets$ \Call{Smoother}{$\M u_l$, $\M f_l$, $\ns[2,l]$}
         \Comment{Post-smoothing}
   \EndFor
   \State \textbf{return} $\M u \gets \M u_L$
\EndFunction
\end{algorithmic}
\end{algorithm}

\subsection{Preconditioned conjugate gradients}
\label{sec:mgcg}

For enhancing robustness and efficiency multigrid methods can be accelerated by Krylov subspace methods \cite{TOS00}.
In the present case, with symmetric system matrices on all grid levels, one would favor preconditioned conjugate gradients.
Unfortunately, weighted additive Schwarz and multiplicative Schwarz with uneven iteration count are both non-symmetric and hence affect the symmetry of MG as well.
As a remedy, it is possible to symmetrize the weighting method or to use GMRES instead of CG for acceleration.
According to \citet{LNS08}, however, symmetrization can deteriorate the efficiency of the method.
This detrimental behavior was confirmed in own tests and, hence, the symmetrized method is not considered here.

Recently, generalizations of the conjugate gradient method have been developed that allow for relaxing some restrictions of standard CG and, thus, promise a cheaper alternative to GMRES.
The use of inaccurately solved and non-symmetric preconditioners
in CG-like methods has been justified, e.g., in \cite{GY99,Not00,Bla02}.
Moreover, \citet{BDN15} demonstrated the suitability of the so-called flexible PCG in conjunction with non-symmetric multigrid preconditioners.
Following this approach, we use the MGCG method summarized in Algorithm~\ref{alg:MGCG}.
This method is equivalent to the flexible PCG of \citet{Not00}, but can be regarded also as a variant of the inexact PCG proposed by \citet{GY99}.
The main difference to standard PCG consists in the application of the Polak-Rib\`iere formula for $\beta$, instead of the Fletcher-Reeves formula, on line~\ref{alg:MGCG:beta} of the algorithm.
We also note that, as before with the coarse problem, the right side $\M f$ must be in the null space of $\M A$ if the system is singular.

\begin{algorithm}[t]
\caption{Inexact multigrid preconditioned conjugate gradients.}
\label{alg:MGCG}
\begin{algorithmic}[1]
\Function{MGCG}{$\M u$, $\M f$, $\M\ns$, $i_{\max}$, $r_{\max}$}
   \State $\M r_{\text{old}} \gets \MS 0$
   \State $\M r \gets \M f - \M A \M u$
   \State $\M p \gets$ \Call{MultigridCycle}{$\MS 0$, $\M r$, $\M\ns$}
   \State $\delta \gets \transpose{\M p} \M r$
   \For{$i = 1, i_{\max}$}
      \State $\M q \gets \M A \M p$
      \State $\alpha \gets \delta / (\transpose{\M p} \M q)$
      \State $\M u \gets \M u + \alpha \M p$
      \State $\M r \gets \M r - \alpha \M q$
      \State \textbf{if} $\lVert \M r \rVert \le r_{\max}$ \textbf{exit}
      \State $\M z \gets$ \Call{MultigridCycle}{$\MS 0$, $\M r$, $\M\ns$}
      \State $\beta \gets \transpose{\M q}(\M r - \M r_{\text{old}}) / \delta$ %
        \label{alg:MGCG:beta}
      \State $\M p \gets \M z + \beta \M p$
      \State $\delta \gets \transpose{\M z} \M r$
      \State $\M r_{\text{old}} \gets \M r$
   \EndFor
   \State \textbf{return} $\M u$
\EndFunction
\end{algorithmic}
\end{algorithm}


\section{Results}
\label{sec:results}

Numerical tests were performed to assess the influence of weighting, overlap and cycling strategy on the computational efficiency and robustness of MG and MGCG.
The methods were implemented in Fortran and compiled using the GNU compiler collection 6.0 with -O3.

All results are based on test cases with the source $f$ evaluated analytically from the exact solution and starting from a random initial guess confined to the interval \mbox{[0,1]}.
The convergence speed is evaluated using
the number of cycles $n_{10}$ needed to reduce the norm of the residual by a factor of $10^{10}$ and the average logarithmic convergence rate according to \citet{Var00}
\begin{equation}
\bar r
= \frac{1}{n} \log_{10} \frac{\lVert \M r^{(0)} \rVert}{\lVert \M r^{(n)} \rVert}
\, ,
\end{equation}
where $\M r^{(n)}$ is the Euclidean norm of the residual vector after the
$n$th cycle.
Note that $n_{10}$ is nearest integer greater than or equal to $10/\bar r$.

As an efficiency measure we define the approximate number of operator applications required for reducing the residual by a factor of $10^k$,
\begin{equation}
\bar\omega_k = \frac{k}{\bar r} \frac{W_{\text{cyc}}}{W_{\text{op}}}
\, ,
\end{equation}
where
$W_{\text{cyc}}$ is the cost for one cycle and
$W_{\text{op}}$ the cost for one application of the system matrix $\M A$.
Exploiting sum factorization \cite{DFM02,KS05}, $W_{\text{op}}$ can be estimated as
$2 \np^3 \nel$, where \mbox{$\np=p+1$}
and $\nel$ is the number of elements.
According to Sec.~\ref{sec:Schwarz}, the cost of one Schwarz iteration is approximately \mbox{$2(\np + 2\no)^3\nel$}.
Assuming a maximum relative overlap of $\no/\np$ this yields the estimate
\begin{equation}
W_{\text{cyc}}
= \left[ 4 \Big(1 + 2\frac{\no}{\np}\Big)^3 \cs\ns
       + 2 \cs + c_{\textsc{cg}}
       \right] \, \np^3 \nel
\, ,
\end{equation}
where
$\ns$ is the number of pre- and post-smoothing steps on the finest level,
\mbox{$\cs = 4/3$} for the classical V-cycle and
\mbox{$\cs = 2$} for a variable V-cycle doubling the number of smoothing steps
when changing to the next lower level, and
\mbox{$c_{\textsc{cg}} = 2$} is the extra cost for conjugate gradients with MGCG.
Since the bracketed term is constant, the overall cost of one multigrid cycle scales approximately with $p N$,
where ${N = p^2 \nel}$ denotes the number of unknowns.
Among the multigrid components, the smoother is by far the most expensive part, accounting for about 80 to 90 percent of the total cost in typical applications.

\subsection{Weighting and overlap}
\label{sec:results:weighting}

We consider the Poisson problem \eqref{eq:poisson} in the domain ${\Omega=[0,2]^2}$, which is uniformly subdivided in \mbox{$8\times8$} square elements with order $p$ ranging from $4$ to $32$.
The right hand side is chosen to match the exact solution ${u=\sin(\pi x)\sin(\pi y)}$.
In the first test series we set the overlap to \mbox{$\no=1$} on all levels \mbox{$l>0$}.
Table~\ref{tab:results:weighting:no=1} shows the measured convergence rates for MG with one pre-smoothing.
Column "$w_a$" corresponds to the weighted additive Schwarz method using the arithmetic mean in overlap areas.
Compared to \cite{LF05} our results agree well for \mbox{$p=4$}, but show a faster convergence with higher polynomial orders.
This could be attributed to using periodic instead of Dirichlet boundary conditions.

The remaining columns in Tab.~\ref{tab:results:weighting:no=1} display the convergence rates for  additive Schwarz smoothing with the gradual weighting introduced in Sec.~\ref{sec:Schwarz} and multiplicative Schwarz.
In comparison to the arithmetic mean, weighting using a smooth -- cubic, quintic or 7\textsuperscript{th} order -- shape function roughly doubled the convergence rate for orders 4, 8 and 16, while linear and top hat weighting yielded a lower, but still remarkable improvement.
As expected, the multiplicative Schwarz smoother attained the fastest convergence.
At \mbox{$p=32$} all methods suffer a serious performance degradation, except for arithmetically weighted Schwarz, which nonetheless remains the slowest.

\begin{table}
\caption{Convergence rates $\bar r$ for MG with additive Schwarz smoothing using one pre- and no post-smoothing steps, $\nel=8\times8$ elements and a fixed overlap of $\no=1$.
The weighting methods are referred to as defined in Tab.~\ref{tab:weighting}.
Results for the multiplicative smoother (mult) are included for comparison.}
\label{tab:results:weighting:no=1}
\centering
\begin{tabular}{rccccccc}
\hline\noalign{\smallskip}
$p$ & \makebox[3em][c]{$w_{a}$}
    & \makebox[3em][c]{$w_1$}
    & \makebox[3em][c]{$w_3$}
    & \makebox[3em][c]{$w_5$}
    & \makebox[3em][c]{$w_7$}
    & \makebox[3em][c]{$w_t$}
    & \makebox[3em][c]{mult}
\\ \noalign{\smallskip}\hline\noalign{\smallskip}
 4  &  0.66  &  0.86  &  1.01  &  1.17  &  1.25  &  0.72  &  1.01  \\
 8  &  0.40  &  0.83  &  1.17  &  1.29  &  1.23  &  0.52  &  1.29  \\
16  &  0.34  &  0.80  &  0.84  &  0.84  &  0.84  &  0.42  &  1.26  \\
32  &  0.32  &  0.43  &  0.43  &  0.43  &  0.43  &  0.38  &  0.76  \\
\noalign{\smallskip}\hline
\end{tabular}
\end{table}

Inspired by these observations, several tests were run with overlaps depending on the polynomial degree on each mesh level.
Table~\ref{tab:results:weighting:no=0.125p} shows the convergence rates for the case
\mbox{$\no[,l]=\lfloor p_l/8 \rfloor$}.
Note that this choice implies \mbox{$\no=0$} for degrees less than 8, while reaching \mbox{$\no=4$} with $p=32$.
As a consequence, the convergence rates for $p=4$ are slightly lower than with \mbox{$\no=1$}, except for multiplicative Schwarz.
For \mbox{$p \ge 16$} the increased overlap yields a considerable speedup.
This improvement is most pronounced for cubic and quintic weighting, which come remarkably close to multiplicative Schwarz.

As a r\'esum\'e of the first study we conclude that
1) gradual weighting with a smooth shape function yields a decisive improvement over arithmetic weighting, and
2) increasing the overlap with growing $p$ is crucial for robustness.

\begin{table}
\caption{Convergence rates for a level-dependent overlap of
$\no[,l]=\lfloor p_l/8 \rfloor$.
For caption see Tab~\ref{tab:results:weighting:no=1}.}
\label{tab:results:weighting:no=0.125p}
\centering
\begin{tabular}{rccccccc}
\hline\noalign{\smallskip}
$p$ & \makebox[3em][c]{$w_{a}$}
    & \makebox[3em][c]{$w_1$}
    & \makebox[3em][c]{$w_3$}
    & \makebox[3em][c]{$w_5$}
    & \makebox[3em][c]{$w_7$}
    & \makebox[3em][c]{$w_t$}
    & \makebox[3em][c]{mult}
\\ \noalign{\smallskip}\hline\noalign{\smallskip}
 4  &  0.63  &  0.91  &  0.98  &  0.96  &  0.79  &  0.31  &  1.03  \\
 8  &  0.40  &  0.75  &  1.06  &  1.28  &  1.28  &  0.64  &  1.30  \\
16  &  0.51  &  1.07  &  1.36  &  1.28  &  1.12  &  0.53  &  1.40  \\
32  &  0.71  &  1.39  &  1.48  &  1.50  &  1.51  &  0.19  &  1.56  \\
\noalign{\smallskip}\hline
\end{tabular}
\end{table}

\subsection{Robustness and efficiency}
\label{sec:robustness}

Next we investigate robustness with respect to the mesh size and aspect ratio.
First, MG with one pre-smoothing is applied on uniform meshes consisting of $4^2$ to $1024^2$ elements with $p$ ranging from 4 to 32 and up to four million unknowns.
Table~\ref{tab:robustness:problem size} compiles the results for quintically weighted and multiplicative Schwarz smoothers with overlap ${\no[,l]=\lceil p_l/8 \rceil}$.
Except in coarse quadrangulations, where periodicity can induce interference effects, the convergence characteristics are virtually independent of the number of elements $\nel$.
The convergence rate $\bar r$ shows a moderate growth for increasing order $p$ and is similar for both smoothers, with a slight advantage for the weighted additive Schwarz method.
As a consequence, the equivalent number of operator applications required for reducing the residual by an order of magnitude drops almost to one third when increasing $p$ from 4 to 32 and, thus, mitigates the higher operator cost per DOF.

In the second test we fixed the mesh to \mbox{$16\times16$} elements of order \mbox{$p=16$}, but  increased the aspect ratio \mbox{$AR=\Delta x / \Delta y$} by enlarging the domain into the $x$ direction, i.e., ${\Omega=[0,2AR]\times[0,2]}$.
Table~\ref{tab:robustness:aspect ratio} reports the results for MG and MGCG using additive weighted Schwarz with $w_5$, ${\no[,l]=\lceil p_l/8 \rceil}$ and one pre-smoothing step.
As expected, the stand-alone MG performs well for small aspect ratios, but degrades for \mbox{$AR>2$}.
MGCG is slightly less efficient than MG for \mbox{$AR\le2$}, but proves more robust at higher aspect ratios.
At \mbox{$AR=8$} it converges approximately twice as fast as MG.

\begin{table}[p]
\caption{Robustness with respect to problem size: MG using additive and multiplicative Schwarz smoothers with overlap ${\no[,l]=\lceil p_l/8 \rceil}$.
}
\label{tab:robustness:problem size}
\centering
\begin{tabular}{rrcccccccc}
\hline\noalign{\smallskip}
     &   && \multicolumn{3}{c}{\makebox[9em][c]{MG(1,0), add, $w_5$}}
         && \multicolumn{3}{c}{\makebox[9em][c]{MG(1,0), mult}}
\\ \noalign{\smallskip}\cline{4-6} \cline{8-10}\noalign{\smallskip}
 $p$ & $\sqrt{\nel}$ && $\bar r$ & $n_{10}$ & $\bar\omega_1$
                     && $\bar r$ & $n_{10}$ & $\bar\omega_1$ \\
\noalign{\smallskip}\hline\noalign{\smallskip}
  4  &   $32$  &&  1.17  &   9  &  9.2  &&  0.87  &  12  & 12.4  \\
     &   $64$  &&  1.17  &   9  &  9.3  &&  0.86  &  12  & 12.6  \\
     &  $128$  &&  1.17  &   9  &  9.3  &&  0.85  &  12  & 12.7  \\
     &  $256$  &&  1.17  &   9  &  9.3  &&  0.85  &  12  & 12.7  \\
\noalign{\smallskip}\hline\noalign{\smallskip}
  8  &   $16$  &&  1.30  &   8  &  5.4  &&  1.28  &   8  &  5.5  \\
     &   $32$  &&  1.29  &   8  &  5.4  &&  1.26  &   8  &  5.5  \\
     &   $64$  &&  1.29  &   8  &  5.4  &&  1.26  &   8  &  5.5  \\
     &  $128$  &&  1.28  &   8  &  5.4  &&  1.26  &   8  &  5.5  \\
\noalign{\smallskip}\hline\noalign{\smallskip}
 16  &    $8$  &&  1.33  &   8  &  5.1  &&  1.44  &   7  &  4.7  \\
     &   $16$  &&  1.37  &   8  &  4.9  &&  1.42  &   8  &  4.8  \\
     &   $32$  &&  1.36  &   8  &  5.0  &&  1.46  &   7  &  4.6  \\
     &   $64$  &&  1.36  &   8  &  5.0  &&  1.46  &   7  &  4.6  \\
\noalign{\smallskip}\hline\noalign{\smallskip}
 32  &    $4$  &&  1.90  &   6  &  3.5  &&  1.65  &   7  &  4.0  \\
     &    $8$  &&  1.58  &   7  &  4.2  &&  1.59  &   7  &  4.2  \\
     &   $16$  &&  1.87  &   6  &  3.6  &&  1.63  &   7  &  4.1  \\
     &   $32$  &&  1.93  &   6  &  3.4  &&  1.64  &   7  &  4.0  \\
     &   $64$  &&  1.93  &   6  &  3.4  &&  1.65  &   7  &  4.0  \\
\noalign{\smallskip}\hline
\end{tabular}
\end{table}

\begin{table}
\caption{Robustness with respect to aspect ratio: MG versus MGCG using additive Schwarz with $w_5$ and overlap ${\no[,l]=\lceil p_l/8 \rceil}$.
}
\label{tab:robustness:aspect ratio}
\centering
\begin{tabular}{rccccccccc}
\hline\noalign{\smallskip}
     &   && \multicolumn{3}{c}{\makebox[9em][c]{MG(1,0), add, $w_5$}}
         && \multicolumn{3}{c}{\makebox[9em][c]{MGCG(1,0), add, $w_5$}}
\\ \noalign{\smallskip}\cline{4-6}\cline{8-10}\noalign{\smallskip}
 $p$ & $AR$ && $\bar r$ & $n_{10}$ & $\bar\omega_1$
            && $\bar r$ & $n_{10}$ & $\bar\omega_1$ \\
\noalign{\smallskip}\hline\noalign{\smallskip}
  4  &  1  &&  1.17  &   9  &   9.2  &&  1.30  &   8  &   9.3  \\
     &  2  &&  0.99  &  11  &  11.0  &&  1.10  &  10  &  11.0  \\
     &  4  &&  0.39  &  26  &  27.8  &&  0.59  &  17  &  20.5  \\
     &  8  &&  0.12  &  85  &  91.2  &&  0.28  &  36  &  42.5  \\
\noalign{\smallskip}\hline\noalign{\smallskip}
  8  &  1  &&  1.30  &   8  &   5.4  &&  1.33  &   8  &   6.1  \\
     &  2  &&  0.86  &  12  &   8.2  &&  1.03  &   8  &   7.9  \\
     &  4  &&  0.43  &  24  &  16.3  &&  0.65  &  16  &  12.5  \\
     &  8  &&  0.16  &  63  &  43.9  &&  0.34  &  30  &  23.8  \\
\noalign{\smallskip}\hline\noalign{\smallskip}
 16  &  1  &&  1.37  &   8  &   4.9  &&  1.55  &   7  &   5.1  \\
     &  2  &&  0.95  &  11  &   7.1  &&  1.14  &   9  &   6.9  \\
     &  4  &&  0.50  &  20  &  13.5  &&  0.72  &  14  &  10.8  \\
     &  8  &&  0.17  &  59  &  39.9  &&  0.39  &  26  &  20.1  \\
\noalign{\smallskip}\hline\noalign{\smallskip}
 32  &  1  &&  1.87  &   6  &   3.6  &&  2.01  &   5  &   3.8  \\
     &  2  &&  1.23  &   9  &   5.4  &&  1.42  &   8  &   5.4  \\
     &  4  &&  0.65  &  16  &  10.3  &&  0.83  &  12  &   9.2  \\
     &  8  &&  0.22  &  46  &  30.4  &&  0.44  &  23  &  17.5  \\
\noalign{\smallskip}\hline
\end{tabular}
\end{table}

\begin{figure}
\centering
\subfloat[Average logarithmic convergence rate]
  {\hspace*{2em}\includegraphics[width=0.65\textwidth]{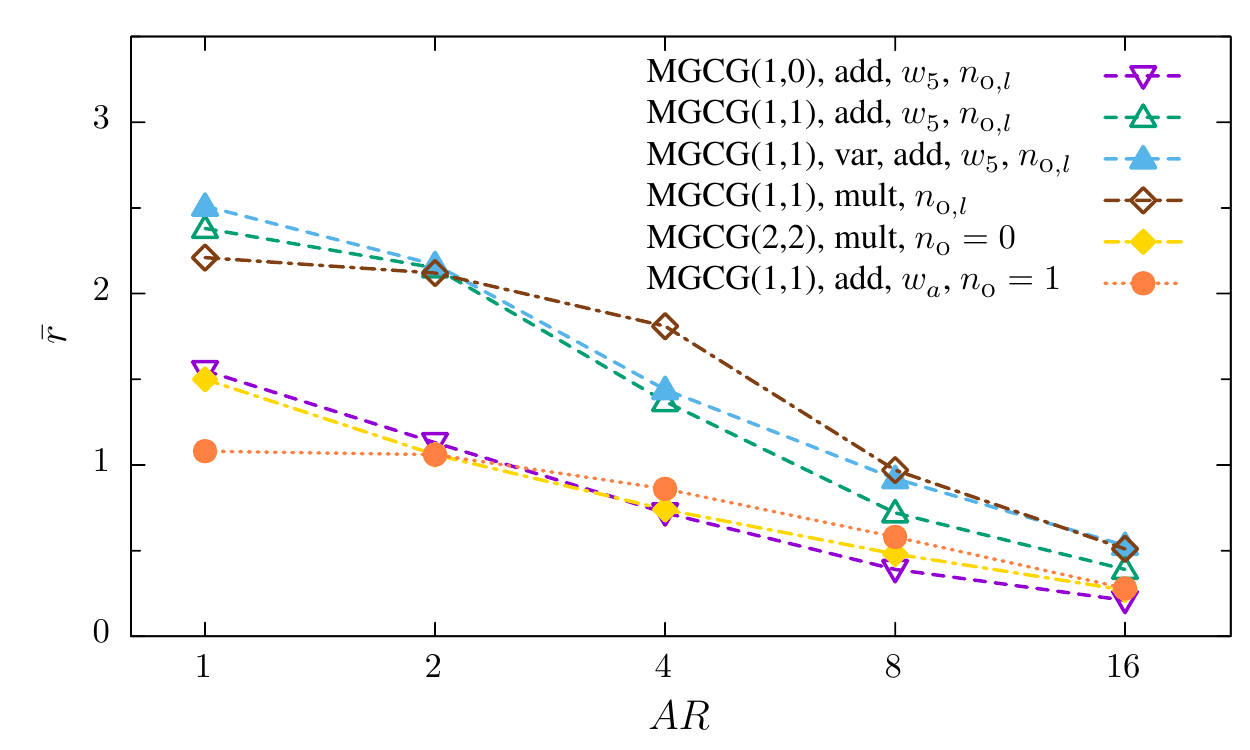}\hspace*{2em}
   \label{fig:p16:ar:r}} \\
\subfloat[Operator applications required for 10\textsuperscript{1} residual reduction]
  {\hspace*{2em}\includegraphics[width=0.65\textwidth]{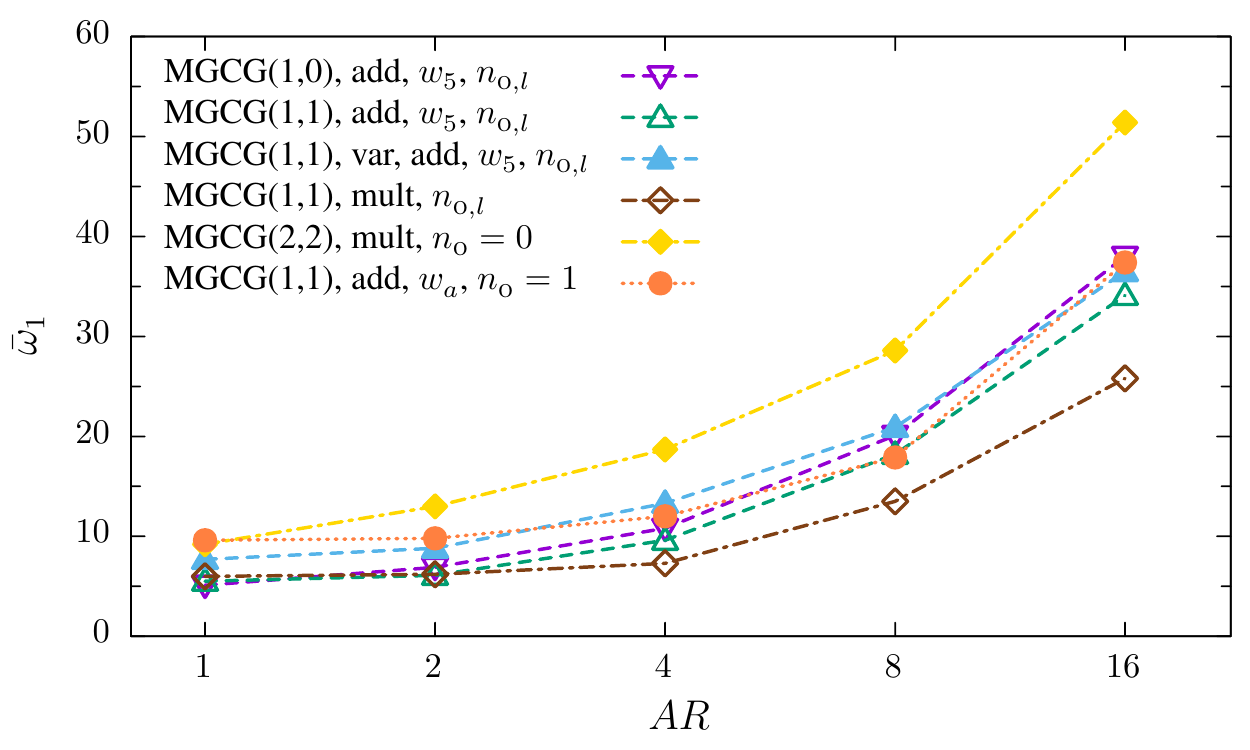}\hspace*{2em}
   \label{fig:p16:ar:omega}} \\
\subfloat[Solver runtime for 10\textsuperscript{10} residual reduction]
  {\hspace*{2em}\includegraphics[width=0.65\textwidth]{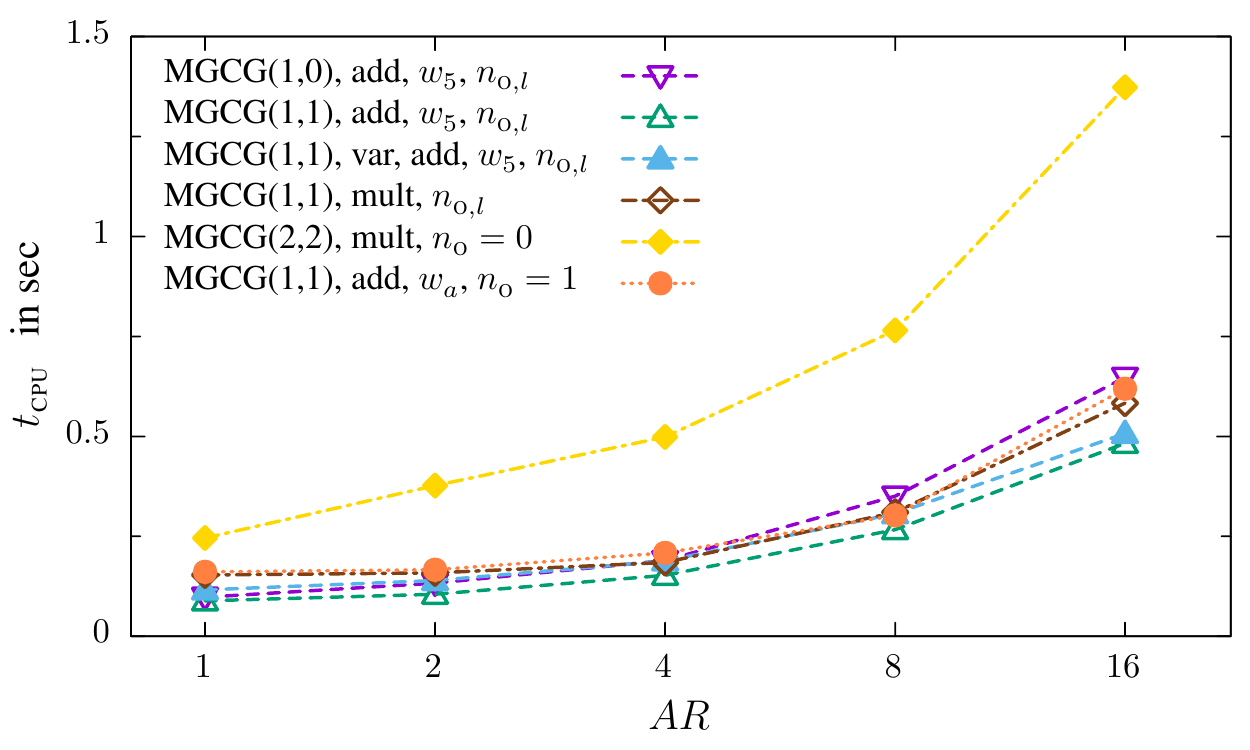}\hspace*{2em}
   \label{fig:p16:ar:t}}
\caption{Performance of MGCG for different aspect ratios. All cases use $16\times16$ elements of order $p=16$ and ${\no[,l]=\lceil p_l/8 \rceil}$, if not specified otherwise.}
\label{fig:p16:ar}
\end{figure}

While these observations hold almost uniformly for all orders $p$ considered, it remains to investigate the impact of solver parameters such as smoothing steps and overlap.
Figure~\ref{fig:p16:ar} presents selected results of the corresponding study for \mbox{$p=16$} and aspect ratios \mbox{$AR=1$} to $16$.
In particular we considered several variants of MGCG(1,1), each applying one pre- and one post-smoothing. In one case, indicated by "var", we employed a variable V-cycle in which the number of smoothing steps doubles with each coarser level, i.e. \mbox{$\ns[1,l]=\ns[2,l]=2^{L-l}$}.
The study included quintically weighted additive ("add, $w_5$") as well as multiplicative ("mult") Schwarz smoothers with a level-dependent overlap of ${\no[,l]=\lceil p_l/8 \rceil}$.
Additionally we tested multiplicative Schwarz with \mbox{$\no=0$} and \mbox{$\ns[1]=\ns[2]=2$}, which corresponds to the method of \citet{JK11},
and the arithmetically averaged additive Schwarz smoother using a constant overlap of \mbox{$\no=1$}.
Figure~\ref{fig:p16:ar:r} depicts the achieved convergence rates.
Compared to the case of only one smoothing, the additional post-smoothing raises $\bar r$ by a factor between 1.5 and 2, which is well in the expected range.
Switching to multiplicative Schwarz yields an even higher gain for increasing aspect ratios.
A similar effect is achieved using additive Schwarz with the variable V-cycle.
MGCG(2,2) with zero overlap attains a convergence rate similar to MGCG(1,0) with level-dependent overlap.
The arithmetically averaged Schwarz method with two smoothing steps falls about two thirds behind the quintically weighted method with only one smoothing for \mbox{$AR=1$}, but gains a slight advantage over the latter for higher aspect ratios.

As the convergence rate does not account for the cost, it is of limited value when comparing methods of different computational complexity. A better measure is the equivalent number of operator applications required for reducing the residual by one order of magnitude, $\bar\omega_1$, which is shown in Fig.~\ref{fig:p16:ar:omega}.
In this metric, the multiplicative MGCG(1,1) with level-dependent overlap performs best, especially for higher aspect ratios.
It is followed by its additive counterpart with quintic weighting, which is at level for \mbox{$AR\le2$}, but needs ca 34 instead of 26 operator applications for \mbox{$AR=16$}.
The comparison also reveals that the benefit of the variable V-cycle is lost due to the higher computational complexity.
Generally, the influence of smoothing and overlap parameters lessens with increasing aspect ratio (exempting the case of \mbox{$\no=0$}), which indicates that the role of the conjugate gradient method gets more important.

Figure~\ref{fig:p16:ar:t} depicts the runtimes measured on a 3.1\,GHz Intel Core i7-5557U  CPU.
Note that MGCG(1,1) with quintic weighting attained the best performance despite its higher operation count in comparison to MGCG(1,1) with multiplicative Schwarz.
This is because the additive Schwarz method evaluates the residual for all elements at once, yielding a single, highly efficient BLAS\,3 operation.
In contrast, multiplicative Schwarz requires a series of local residual updates, which is harder to optimize.
Consistently, the multiplicative MGCG(2,2) with ${\no=0}$ remains the least efficient method for all aspect ratios.
Compared to MGCG(1,1) with ${\no=1}$ and arithmetic weighting, the method with ${\no[,l]=\lceil p_l/8 \rceil}$ and quintic weighting succeeds twice as fast for \mbox{$AR=1$} and still gains 23\% at \mbox{$AR=16$}.
Though other choices may yield even better performance, the study documents that the method is not too sensitive to parameter variations, such that only minor improvements can be expected.

\subsection{Variable diffusion}
Since many problems in physics involve variable coefficients, it is interesting to explore if the multigrid method is capable to retain its efficiency in such applications.
As an example we consider the diffusion equation
\begin{equation}
\label{eq:diffusion}
-\nabla \cdot (\nu \nabla u) = f
\end{equation}
with variable diffusivity $\nu$ in the periodic domain
${\Omega=[0,AR]\times[0,1]}$.
From a physical perspective it seems reasonable that the solution and the diffusivity vary on similar scales.
Following this idea we set ${u = \sin(2\pi x) \sin(2\pi y)}$ and
\begin{equation}
\label{eq:diffusion:nu}
\nu = 1 + \hat\nu \sin\big(2\pi (x - s)\big) \sin \big(2\pi (y - s)\big),
\end{equation}
where
$\hat\nu$ is the amplitude and
$s$ the shift of the diffusivity fluctuation.
According to preliminary studies, a non-zero shift poses an additional difficulty to the solver.
Taking this into account, ${s = 0.2}$ is chosen in all tests reported below.
The source is analytically computed from
${f = -(\nu \nabla^2 u + \nabla \nu \cdot \nabla u)}$.

Discretization using rectangular spectral elements yields the linear system
\begin{equation}
\label{eq:diffusion:discrete system}
\M B \M u = \M f
\, ,
\end{equation}
where $\M B(\nu)$ represents the discrete diffusion operator or, equivalently,
\begin{equation}
\label{eq:diffusion:residual}
\M B \Delta \M u = \M f - \M B \M{\tilde u} = \M{\tilde r}
\, ,
\end{equation}
for the correction $\Delta \M u$ to a given approximation $\M{\tilde u}$.
Application of the Schwarz method described in Sec.~\ref{sec:Schwarz} leads to the local correction equation
\begin{equation}
\label{eq:subproblem:diffusion}
\M B_{ss} \Delta \M u_s = \M r_s \, ,
\end{equation}
where $\M B_{ss}$ and $\M r_s$ are the diffusion operator and, respectively, the residual
restricted to the subdomain $\Omega_s$.
In comparison to the subdomain problem for Poisson case \eqref{eq:subproblem}, Equation \eqref{eq:subproblem:diffusion} is more expensive to solve, because the fast diagonalization technique is no longer applicable.
Yet, the smoothing property is more important for multigrid than accurate solution of the subproblems.
This motivates the reintroduction of the discrete Laplacian by approximating the restricted diffusion operator on the left side of \eqref{eq:subproblem:diffusion} by
${\M B_{ss} \approx \bar\nu_{s}} \M A_{ss}$, where the diffusivity $\bar\nu_{s}$ is assumed to be constant in $\Omega_s$.
For simplicity, $\bar\nu_{s}$ is set to the average of $\nu$ over the embedded element.
The correction can then be approximated as
\begin{equation}
\label{eq:subproblem:solution}
\Delta \M u_s \approx \frac{1}{\bar\nu_{s}}\M A_{ss}^{-1} \M r_s \, ,
\end{equation}
where, again, $\M A_{ss}^{-1}$ stands for the application of the factored inverse obtained from fast diagonalization.
As a result, the solution techniques developed in Sec.~\ref{sec:methods} can be utilized with no change except for the residual evaluation.

\begin{figure}[t!]
\centering
\includegraphics[width=0.65\textwidth]{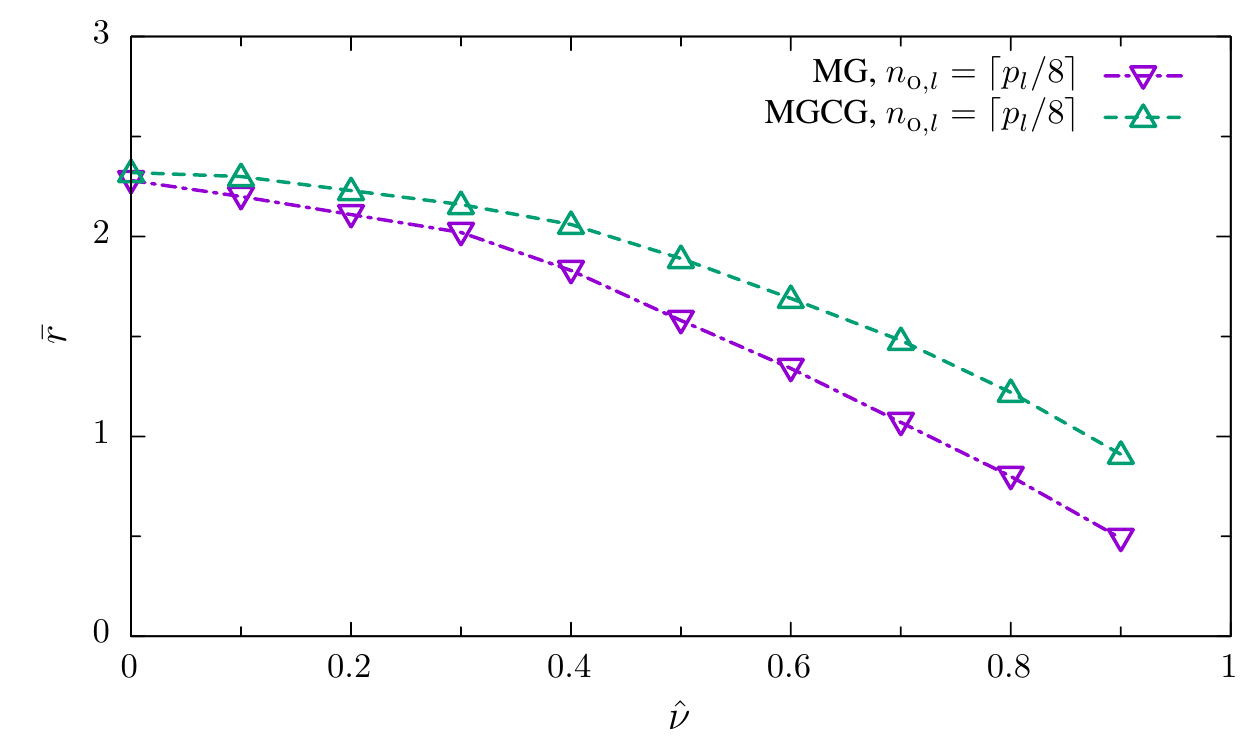}
\caption{%
MG and MGCG convergence rates for different diffusivity fluctuation amplitudes.
Discretization is based on an isotropic mesh comprising $8^2$ elements of order ${p=16}$.
One pre- and post-smoothing with an overlap of ${\no[,l]=\lceil p_l/8 \rceil}$ and quintic weighting were applied in both cases.
\label{fig:diffusion:cv}}
\end{figure}

The performance of the scheme was studied in two test series.
In the first series, the aspect ratio was fixed to ${AR=1}$ and the domain  ${\Omega = [0,1]^2}$ decomposed into $8^2$ square elements of order ${p = 16}$.
The diffusivity fluctuation amplitude $\hat\nu$ was gradually increased from 0 to 0.9, where the latter corresponds to variations of the magnitude up to 90 percent.
Figure~\ref{fig:diffusion:cv} shows the measured convergence rates for MG und MGCG using one pre- and post-smoothing based on additive Schwarz with a level-dependent overlap of ${\no[,l]=\lceil p_l/8 \rceil}$ and quintic weighting.
The results indicate that MG retains its efficiency up to fluctuation amplitudes of about 30 percent, but then degrades with rising $\hat\nu$.
As expected, Krylov acceleration improves the robustness, such that MGCG achieves ${\bar r = 0.91}$ for ${\hat\nu = 0.9}$, which is nearly twice the convergence rate obtained with MG.
Compared to ${\hat\nu = 0}$, this corresponds to an increase of the cycle count and, hence, in runtime, by a factor of just $2.2$.
A similar behavior was observed for polynomial orders ${p=4}$, $8$ and $32$.

\begin{figure}
\centering
\includegraphics[width=0.65\textwidth]{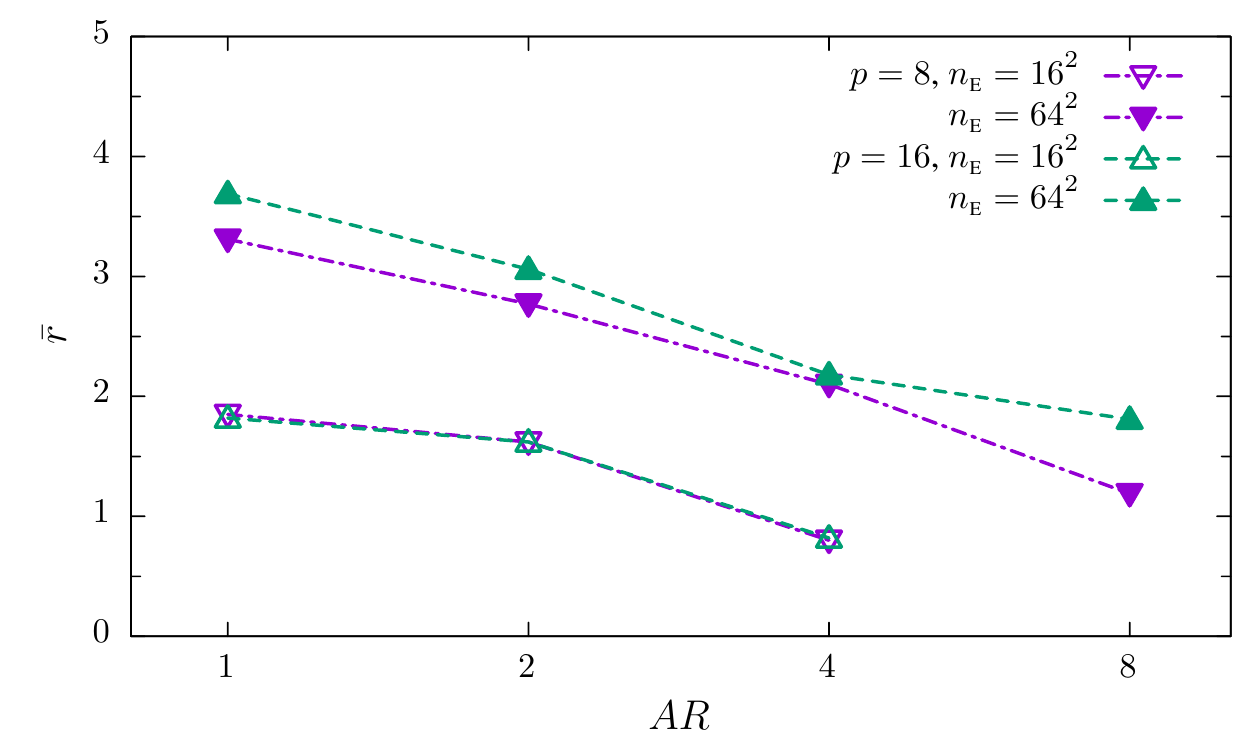}
\caption{%
MGCG convergence rate for variable diffusivity on anisotropic meshes with increasing aspect ratio using a variable V-cycle with one pre- and post-smoothing, overlap ${\no[,l]=\lceil p_l/2 \rceil}$ and quintic weighting.}
\label{fig:diffusion:ar}
\end{figure}

In the second test series, we increased the domain extension in the $x$-direction,
while keeping the diffusivity fluctuation amplitude $\hat\nu$ at a constant level.
The number of elements is fixed and identical in both directions, such that the element aspect ratio equals $AR$.
For achieving a robustness similar to the Poisson case it proved necessary to increase the subdomain overlap with growing $\hat\nu$.
Figure~\ref{fig:diffusion:ar} shows the convergence rates $\bar r$ for MGCG using a variable V-cycle and additive Schwarz smoothing for an amplitude of 90 percent, which represents the most challenging test in the series.
Comparing the results for ${p=8}$ and ${p=16}$ one observes that $\bar r$ strongly depends on the quadrangulation, but only marginally on the polynomial order.
Using a finer mesh yields considerably higher convergence rates and better robustness.
Orders $4$ and ${32}$ fit nicely into this picture, but are not shown for clarity.
The congruence of different orders using the same mesh suggests, that the performance depends on how well the diffusivity fluctuation is resolved by the element mean values adopted for $\bar\nu_s$.
This presents a possible limitation of the approach, which needs further consideration in subsequent work.
Nevertheless, the study demonstrates the suitability of the proposed method for problems involving variable diffusivity, as long as the latter is sufficiently resolved.

\section{Conclusions}
\label{sec:conclusions}

We have developed a nonuniformly weighted additive Schwarz method acting as the smoother in multigrid solvers for the spectral element discretization of the Poisson equation in $\mathbb R^2$.
The method generalizes the Schwarz/multigrid method proposed in \cite{LF05} and was inspired from weighting techniques devised in \cite{HSN13}.
In each step, it determines the solution for a subdomain corresponding to an extended element region.
These local solutions are blended according to a polynomial shape function which features a smooth transition from zero at the border toward one in the core of the subdomain.
As an alternative we considered a multiplicative Schwarz method with no weighting required.
Both Schwarz methods were integrated in a polynomial multigrid method which, in turn, was embedded in a preconditioned CG method.

The performance of these methods was assessed in a series of numerical experiments with ansatz orders $p$ ranging from 4 to 32 and up to \mbox{$\nel=256^2$} elements of aspect ratios $AR$ from 1 to 16.
For unit-aspect ratio elements the proposed weighting improved the logarithmic MG convergence rate and reduced the cost by a factor of 1.5 to 3 in comparison to the original method.
The study indicates that for robustness the subdomain overlap has to be bounded, i.e., the number $\no$ of node layers adopted from neighbor elements must grow with increasing order.
Thus, with MG, the number of layers varies from level to level.
A reasonable choice is to use an overlap of ${\lceil p_l/8\rceil}$ layers, where $p_l$ denotes the polynomial order on level $l$.
The resulting multigrid method is robust with respect to the mesh size, i.e. $p$ and $\nel$, but degrades with increasing aspect ratio.
This behavior can be mitigated by Krylov subspace acceleration:
Using MG as a preconditioner for the inexact conjugate gradient method \cite{GY99} improves the convergence rate for higher aspect ratios considerably.

Finally, it has been shown that the proposed multigrid method is easily adapted and well suited for solving diffusion problems with varying coefficients, provided the mesh is fine enough to approximate diffusivity fluctuations by element mean values.
Improving the treatment of variable coefficients and extending the approach to three space dimensions are topics of ongoing and future work.

\section*{Acknowledgements}
Funding by German Research Foundation (DFG) in frame of the project STI 157/4-1 is gratefully acknowledged.



\begin{thebibliography}{29}
\providecommand{\natexlab}[1]{#1}
\providecommand{\url}[1]{\texttt{#1}}
\expandafter\ifx\csname urlstyle\endcsname\relax
  \providecommand{\doi}[1]{doi: #1}\else
  \providecommand{\doi}{doi: \begingroup \urlstyle{rm}\Url}\fi

\bibitem[Bastian et~al.(2012)Bastian, Blatt, and Scheichl]{BBS12}
P.~Bastian, M.~Blatt, and R.~Scheichl.
\newblock Algebraic multigrid for discontinuous {G}alerkin discretizations of
  heterogeneous elliptic problems.
\newblock \emph{Numer. Linear Algebra Appl.}, 19\penalty0 (2):\penalty0
  367--388, 2012.

\bibitem[Blaheta(2002)]{Bla02}
R.~Blaheta.
\newblock {GPCG}--generalized preconditioned {CG} method and its use with
  non-linear and non-symmetric displacement decomposition preconditioners.
\newblock \emph{Numerical Linear Algebra with Applications}, 9\penalty0
  (6-7):\penalty0 527--550, 2002.

\bibitem[Bouwmeester et~al.(2015)Bouwmeester, Dougherty, and Knyazev]{BDN15}
H.~Bouwmeester, A.~Dougherty, and A.~V. Knyazev.
\newblock Nonsymmetric preconditioning for conjugate gradient and steepest
  descent methods.
\newblock \emph{Procedia Computer Science}, 51:\penalty0 276--285, 2015.
\newblock ISSN 1877-0509.

\bibitem[Bramble(1995)]{Bra93}
J.~Bramble.
\newblock \emph{Multigrid methods}.
\newblock Pitman Res. Notes Math. Ser. 294. Longman Scientific \& Technical,
  Harlow, UK, 1995.

\bibitem[Brenner and Zhao(2005)]{BZ05}
S.~C. Brenner and J.~Zhao.
\newblock Convergence of multigrid algorithms for interior penalty methods.
\newblock \emph{Appl. Num. Anal. Comp. Math.}, 2\penalty0 (1):\penalty0 3--18,
  2005.

\bibitem[Couzy and Deville(1995)]{CD95}
W.~Couzy and M.~O. Deville.
\newblock A fast {S}chur complement method for the spectral element
  discretization of the incompressible {N}avier-{S}tokes equations.
\newblock \emph{J. Comput. Phys.}, 116:\penalty0 135--142, January 1995.

\bibitem[Deville et~al.(2002)Deville, Fischer, and Mund]{DFM02}
M.~O. Deville, P.~F. Fischer, and E.~H. Mund.
\newblock \emph{High-Order Methods for Incompressible Fluid Flow}, volume~1.
\newblock Cambridge University Press, 2002.

\bibitem[Fischer and Lottes(2004)]{FL04}
P.~F. Fischer and J.~W. Lottes.
\newblock Hybrid {S}chwarz-multigrid methods for the spectral element method:
  Extensions to {N}avier-{Stokes}.
\newblock In \emph{Domain Decomposition Methods in Science and Engineering
  Series}, pages 35--49. Springer, 2004.

\bibitem[Golub and Ye(1999)]{GY99}
G.~H. Golub and Q.~Ye.
\newblock Inexact preconditioned conjugate gradient method with inner-outer
  iteration.
\newblock \emph{SIAM Journal on Scientific Computing}, 21\penalty0
  (4):\penalty0 1305--1320, Dec. 1999.
\newblock ISSN 10648275.

\bibitem[Guermond et~al.(2006)Guermond, Minev, and Shen]{Guermond2006}
J.~L. Guermond, P.~Minev, and J.~Shen.
\newblock An overview of projection methods for incompressible flows.
\newblock \emph{Comput. Methods Appl. Mech. Eng.}, 195:\penalty0 6011--6045,
  2006.

\bibitem[Haupt et~al.(2013)Haupt, Stiller, and Nagel]{HSN13}
L.~Haupt, J.~Stiller, and W.~Nagel.
\newblock A fast spectral element solver combining static condensation and
  multigrid techniques.
\newblock \emph{J. Comput. Phys.}, 255:\penalty0 384--395, 2013.

\bibitem[Heinrichs(1988)]{Hei88}
W.~Heinrichs.
\newblock Line relaxation for spectral multigrid methods.
\newblock \emph{J. Comput. Phys.}, 77:\penalty0 166--182, 1988.

\bibitem[Janssen and Kanschat(2011)]{JK11}
B.~Janssen and G.~Kanschat.
\newblock Adaptive multilevel methods with local smoothing for \mbox{$H^1$}-
  and \mbox{$H^{\mathrm{curl}}$}-conforming high order finite element methods.
\newblock \emph{{SIAM} J. Sci. Comput.}, 33\penalty0 (4):\penalty0 2095--2114,
  2011.

\bibitem[Kaasschieter(1988)]{Kaa88}
E.~F. Kaasschieter.
\newblock Preconditioned conjugate gradients for solving singular systems.
\newblock \emph{J. Comput. Appl. Math.}, 24:\penalty0 265--275, 1988.

\bibitem[Kanschat(2004)]{Kan04}
G.~Kanschat.
\newblock Multilevel methods for discontinuous {G}alerkin {FEM} on locally
  refined meshes.
\newblock \emph{Computers and Structures}, 82:\penalty0 2437--2445, 2004.

\bibitem[Karniadakis and Sherwin(2005)]{KS05}
G.~E. Karniadakis and S.~J. Sherwin.
\newblock \emph{Spectral/hp Element Methods for Computational Fluid Dynamics}.
\newblock Oxford University Press, 2nd edition, 2005.

\bibitem[Kraus and Tomar(2008)]{KT08}
J.~K. Kraus and S.~K. Tomar.
\newblock A multilevel method for discontinuous {G}alerkin approximation of
  three-dimensional anisotropic elliptic problems.
\newblock \emph{Numer. Linear Algebra Appl.}, 15\penalty0 (5):\penalty0
  417--438, 2008.

\bibitem[Loisel et~al.(2008)Loisel, Nabben, and Szyld]{LNS08}
S.~Loisel, R.~Nabben, and D.~Szyld.
\newblock On hybrid multigrid-{S}chwarz algorithms.
\newblock \emph{J. Sci. Comput.}, 36\penalty0 (2):\penalty0 165--175, 2008.

\bibitem[Lottes and Fischer(2005)]{LF05}
J.~W. Lottes and P.~F. Fischer.
\newblock Hybrid multigrid/{S}chwarz algorithms for the spectral element
  method.
\newblock \emph{J. Sci. Comput.}, 24:\penalty0 45--78, 2005.

\bibitem[Lynch et~al.(1964)Lynch, Rice, and Thomas]{Lyn64}
R.~E. Lynch, J.~R. Rice, and D.~H. Thomas.
\newblock Direct solution of partial difference equations by tensor product
  methods.
\newblock \emph{Numer. Math.}, 6:\penalty0 185--199, 1964.

\bibitem[Maday and Munoz(1988)]{MM88}
Y.~Maday and R.~Munoz.
\newblock Spectral element multigrid. {II}. {T}heoretical justification.
\newblock \emph{J. Sci. Comput.}, 3:\penalty0 323--353, 1988.

\bibitem[Mitchell(2010)]{Mit10}
W.~F. Mitchell.
\newblock The hp-multigrid method applied to hp-adaptive refinement of
  triangular grids.
\newblock \emph{Numer. Linear Algebr.}, 17:\penalty0 211--228, 2010.

\bibitem[Notay(2000)]{Not00}
Y.~Notay.
\newblock Flexible conjugate gradients.
\newblock \emph{SIAM Journal on Scientific Computing}, 22\penalty0
  (4):\penalty0 1444--1460, 2000.

\bibitem[Olson(2007)]{Ols07}
L.~Olson.
\newblock Algebraic multigrid preconditioning of high-order spectral elements
  for elliptic problems on a simplicial mesh.
\newblock \emph{SIAM J. Sci. Comput.}, 29\penalty0 (5):\penalty0 2189--2209,
  2007.

\bibitem[Pasquetti and Rapetti(2009)]{PR09}
R.~Pasquetti and F.~Rapetti.
\newblock p-multigrid method for {F}ekete-{G}auss spectral element
  approximations of elliptic problems.
\newblock \emph{Commun. Comput. Phys.}, 5\penalty0 (5):\penalty0 667--682,
  February 2009.

\bibitem[R{\o}nquist and Patera(1987)]{RP87}
E.~R{\o}nquist and A.~Patera.
\newblock Spectral element multigrid. {I}. {F}ormulation and numerical results.
\newblock \emph{J. Sci. Comput.}, 2:\penalty0 389--406, 1987.

\bibitem[Trottenberg et~al.(2000)Trottenberg, Oosterlee, and
  Sch\"{u}ller]{TOS00}
U.~Trottenberg, C.~W. Oosterlee, and A.~Sch\"{u}ller.
\newblock \emph{Multigrid}.
\newblock Academic Press, 2000.

\bibitem[Varga(2000)]{Var00}
R.~S. Varga.
\newblock \emph{Matrix Iterative Analysis}.
\newblock Springer Ser. Comput. Math. 27. Springer-Verlag, Berlin, 2nd edition,
  2000.

\bibitem[Wang et~al.(2013)Wang, Fidkowski, Abgrall, Bassi, Caraeni, Cary,
  Deconinck, Hartmann, Hillewaert, Huynh, Kroll, May, Persson, van Leer, and
  Visbal]{Wang2013}
Z.~J. Wang, K.~Fidkowski, R.~Abgrall, F.~Bassi, D.~Caraeni, A.~Cary,
  H.~Deconinck, R.~Hartmann, K.~Hillewaert, H.~Huynh, N.~Kroll, G.~May, P.-O.
  Persson, B.~van Leer, and M.~Visbal.
\newblock High-order {CFD} methods: current status and perspective.
\newblock \emph{International Journal for Numerical Methods in Fluids},
  72\penalty0 (8):\penalty0 811--845, 2013.
\newblock ISSN 1097-0363.

\end{thebibliography}
\end{document}